\documentclass[11pt,here]{article}
\usepackage[colorlinks=true,
linkcolor=webgreen,
filecolor=webbrown,
citecolor=webgreen]{hyperref}

\usepackage{amssymb,psfig,epsfig}

\definecolor{webgreen}{rgb}{0,.5,0}
\definecolor{webbrown}{rgb}{.6,0,0}

\usepackage{psfig,epsfig}
\newtheorem{theorem}{Theorem}

\newtheorem{prop}{Proposition}
\newtheorem{exam}{Example}

\def\binom#1#2{{#1}\choose{#2}}

\newcommand{\eqn}[1]{(\ref{#1})}
\newcommand{\hsp}{\hspace*{\parindent}}

\newcommand{\eeq}{\end{equation}}
\newcommand{\beql}[1]{\begin{equation}\label{#1}}
\newcommand{\bsq}{{\vrule height .9ex width .8ex depth -.1ex }}

\newcommand{\sE}{{\cal E}}
\newcommand{\sA}{{\cal A}}
\newcommand{\sM}{{\cal M}}

\newcommand{\ba}{{\bf a}}
\newcommand{\bb}{{\bf b}}
\newcommand{\bm}{{\bf m}}

\newcommand{\dd}{\ldots}

\newcommand{\sB}{{\cal B}}

\makeatletter
\def\@sect#1#2#3#4#5#6[#7]#8{\ifnum #2>\c@secnumdepth
     \def\@svsec{}\else
     \refstepcounter{#1}\edef\@svsec{\csname the#1\endcsname.\hskip .75em }\fi
     \@tempskipa #5\relax
      \ifdim \@tempskipa>\z@
        \begingroup #6\relax
          \@hangfrom{\hskip #3\relax\@svsec}{\interlinepenalty \@M #8\par}%
        \endgroup
       \csname #1mark\endcsname{#7}\addcontentsline
         {toc}{#1}{\ifnum #2>\c@secnumdepth \else
                      \protect\numberline{\csname the#1\endcsname}\fi
                    #7}\else
        \def\@svsechd{#6\hskip #3\@svsec #8\csname #1mark\endcsname
                      {#7}\addcontentsline
                           {toc}{#1}{\ifnum #2>\c@secnumdepth \else
                             \protect\numberline{\csname the#1\endcsname}\fi
                       #7}}\fi
     \@xsect{#5}}
\def\@begintheorem#1#2{\it \trivlist \item[\hskip \labelsep{\bf #1\ #2.}]}
\makeatother

\catcode`\@=11
\renewcommand{\section}{
        \setcounter{equation}{0}
        \@startsection {section}{1}{\z@}{-3.5ex plus -1ex minus
        -.2ex}{2.3ex plus .2ex}{\large\bf}
        }
\catcode`@=12

\begin{document}
\begin{center}
{\Large {\bf A New Operation on Sequences: The Boustrophedon Transform}}\footnote{This paper was published (in a somewhat different form) in {\em J. Combinatorial Theory}, Series A, {\bf 76} (1996), pp. 44--54.}
\\
\vspace{1.5\baselineskip}
{\em J. Millar}\footnote{Present address: Mathematics Department, MIT, Cambridge, MA},
{\em N. J. A. Sloane}\footnote{Present address: AT\&T Shannon Labs, Florham Park, NJ}
and {\em N. E. Young}\footnote{Present address: Akamai Technologies, Cambridge, MA} \\
\vspace*{+.2\baselineskip}
Mathematical Sciences Research Center \\
AT\&T Bell Laboratories \\
Murray Hill, New Jersey 07974 \\
\vspace{1\baselineskip}
January 15, 1996 \\
Enhanced version April 26, 2000 \\
Postscript added June 10, 2002 \\
\vspace{1.5\baselineskip}
{\bf ABSTRACT}
\vspace{.5\baselineskip}
\end{center}
\setlength{\baselineskip}{1.5\baselineskip}

A generalization of the Seidel-Entringer-Arnold method for calculating
the alternating permutation numbers (or secant-tangent numbers)
leads to a new operation on sequences, the boustrophedon transform.
\clearpage
\section{Introduction}
\hsp
Let $E_{n,k}$ $(n \ge k \ge 0)$ denote the number of permutations
of $\{1,2, \dd, n+1 \}$ which alternately fall and rise (always starting with a fall), and start with $k+1$.
These numbers have a long history (see the references), but we follow
Poupard \cite{Pou82} and call them the Entringer numbers.
They satisfy the recurrence \cite[first lemma]{Ent66}
\beql{Eq1}
E_{0,0} =1, ~~E_{n,0} =0 ~ (n \ge 1) , ~~
E_{n+1, k+1} = E_{n+1,k} + E_{n,n-k} ~~ (n \ge k \ge 0) ~.
\eeq
If these numbers are displayed in a triangular array with rows written
alternately right to left and left to right, in boustrophedon (or
``ox-plowing'') manner (sequence
\htmladdnormallink{A8280}{http://www.research.att.com/cgi-bin/access.cgi/as/njas/sequences/eisA.cgi?Anum=008280} in \cite{OEIS}):
\beql{Eq2}
\begin{array}{ccccccccc}
~ & ~ & ~ & ~ & E_{00} \\
~ & ~ & ~ & E_{10} & \to & E_{11} \\
~ & ~ & E_{22} & \leftarrow & E_{21} & \leftarrow & E_{20} \\
~ & E_{30} & \to & E_{31} & \to & E_{32} & \to & E_{33} \\
E_{44} & \leftarrow & E_{43} & \leftarrow & E_{42} & \leftarrow & E_{41} & \leftarrow & E_{40} \\
~ & ~ & ~ & ~ & \cdots
\end{array}
~~=~~
\begin{array}{ccccccccc}
~ & ~ & ~ & ~ & 1 \\
~ & ~ & ~ & 0 & \to & 1 \\
~ & ~ & 1 & \leftarrow & 1 & \leftarrow & 0 \\
~ & 0 & \to & 1 & \to & 2 & \to & 2 \\
5 & \leftarrow & 5 & \leftarrow & 4 & \leftarrow & 2 & \leftarrow & 0 \\
~ & ~ & ~ & ~ & \cdots
\end{array}
\eeq
then the entries are filled in by the rule that each row (after the
zero-th) begins with a 0 and every subsequent entry is the sum of the previous
entry in the same row and the entry above it in the previous row.

The earliest reference we have seen for this elegant observation
is Arnold \cite{Arn91}, who refers to \eqn{Eq2} as the Euler-Bernoulli
triangle, but it may well be of much older origin.
Dumont \cite{Dum95} refers to \eqn{Eq2} as the Seidel-Entringer-Arnold triangle,
referring to Seidel \cite{Sei77}.

The numbers $E_n := E_{n,n}$
(sequence
\htmladdnormallink{A111}{http://www.research.att.com/cgi-bin/access.cgi/as/njas/sequences/eisA.cgi?Anum=000111})
appearing at the ends of the rows in
\eqn{Eq2} give the total number of permutations of $\{1,2, \dd, n\}$ that alternately fall and rise,
i.e. the number of ``down-up permutations'' of $n$ things.
The history of these numbers goes back to Andr\'{e} \cite{And79},
\cite{And81}, \cite{Com74}, \cite{Sch61}.
They have exponential generating function (e.g.f.)
\beql{Eq3}
\sE(x) = \sum_{n=0}^\infty E_n \frac{x^n}{n!} = {\rm sec} \, x + \tan \, x ~.
\eeq
Conway and Guy \cite{CG96} call \eqn{Eq2} the zig-zag triangle and the $E_n$ the zig-zag permutation numbers.
The Entringer numbers have also been shown to enumerate several classes of
rooted planar trees as well as other mathematical objects
\cite{Arn91}, \cite{Arn92}, \cite{Kem33}, \cite{KPP94}, \cite{Pou82}.

Guy \cite{Guy95} observed that if the entries at the beginnings of the rows in \eqn{Eq2} are changed from $1,0,0,0, \dd$ to say $1,1,1,1,1, \dd$, or
$1,2,4,8,16, \dd$, etc.,
then the numbers at the ends of the rows form interesting-looking
sequences not to be found in \cite{SP95}.
(Of course now they are in \cite{OEIS}.)
Using $1,1,1, \dd$ for example the triangle becomes
\beql{Eq4}
\begin{array}{ccccccccccc}
~ & ~ & ~ & ~ & ~ & 1 \\
~ & ~ & ~ & ~ & ~ 1 & ~ & 2 \\
~ & ~ & ~ & 4 & ~ & 3 & ~ & 1 \\
~ & ~ & 1 & ~ & 5 & ~ & 8 & ~ & 9 \\
~ & 24 & ~ & 23 & ~ & 18 & ~ & 10 & ~ & 1 \\
1 & ~ & 25 & ~ & 48 & ~ & 66 & ~ & 76 & ~ & 77 \\
~ & ~ & ~ & ~ & ~ & \cdots
\end{array}
\eeq
yielding the sequence
(\htmladdnormallink{A667}{http://www.research.att.com/cgi-bin/access.cgi/as/njas/sequences/eisA.cgi?Anum=000667})
\beql{Eq15}
1,2,4,9,24,77,294,1309, \dd ~.
\eeq
Guy asked if anything could be said about generating functions or
combinatorial interpretations for these sequences.
The purpose of this note is to answer his question.
\section{The Boustrophedon transform}
\hsp
Given a sequence\footnote{In this paper we consider only integer-valued sequences, although the transformation can be applied to sequences over any ring.}
$\ba = (a_1, a_1, a_2, \dd )$ we define its boustrophedon transform
to be the sequence $\bb = (b_0, b_1, b_2, \dd)$ produced by the
triangle
\beql{Eq9}
\begin{array}{ccccccc}
~ & ~ & ~ & a_0 = b_0 \\ [+.15in]
~ & ~ & a_1 & \to & b_1 = a_0 + a_1 \\ [+.15in]
~ & b_2 = a_1 + a_2 + b_1 & \leftarrow & a_2 + b_1 & \leftarrow & a_2 \\ [+.15in]
a_3 & \to & a_3 + b_2 & \to & a_2 + a_3 + b_1 + b_2 & \to & b_3 = 2a_2 + a_3 + b_1 + b_2 \\ [+.15in]
~ & ~ & ~ & \cdots
\end{array}
\eeq
when it is filled in using the rule described in Section~1.
Formally, the entries $T_{n,k}$ $(n \ge k \ge 0)$ in the triangle
are defined by
\beql{Eq5}
\begin{array}{rll}
T_{n,0} & = & a_n ~~~~(n \ge 0) ~, \\ [+.1in]
T_{n+1,k+1} & = & T_{n+1,k} + T_{n,n-k} ~~~~(n \ge k \ge 0) ~,
\end{array}
\eeq
and then
$$b_n = T_{n,n} ~~~~ (n \ge 0) ~.$$
Although many operations on sequences have been studied in the past (see \cite{BS95} and the references therein), this transformation
appears to have been overlooked.
\begin{theorem}
\label{th1}
The boustrophedon transform $\bb$ of a sequence $\ba$ is given by
\begin{eqnarray}
\label{Eq6}
b_n & = & \sum_{k=0}^n {\binom{n}{k}} a_k E_{n-k} ,~~~~(n \ge 0)~, \\
\label{Eq7}
a_n & = & \sum_{k=0}^n (-1)^{n-k} {\binom{n}{k}} b_k E_{n-k} , ~~~~
(n \ge 0) ~,
\end{eqnarray}
and the e.g.f.'s of $\bb$ and $\ba$ are related by
\beql{Eq8}
\sB (x) = ( {\rm sec} \, x + \tan x) ~ \sA (x) ~.
\eeq
\end{theorem}
\paragraph{Proof.}
We redraw \eqn{Eq9} as a directed graph $\Gamma$ whose nodes are labeled by the numbers $T_{n,k}$ (see Fig.~1).
\begin{figure}[htb]
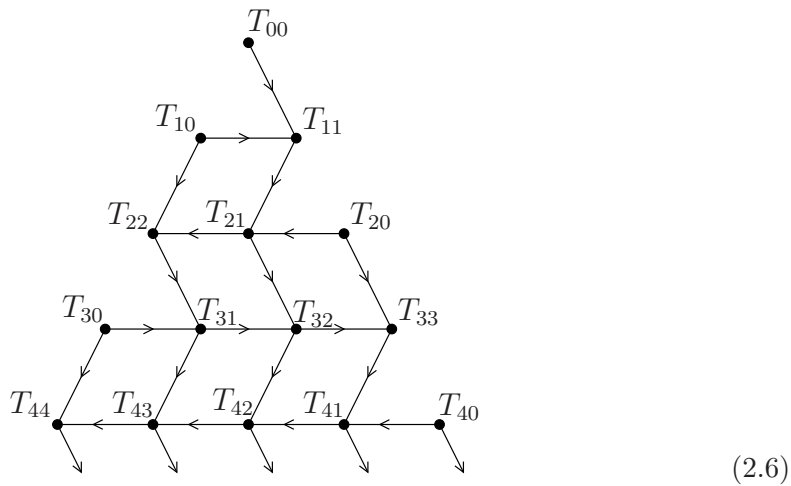

$$
\mbox{
\input bo1.tex
}
\eqno{(2.6)}
$$
\caption{Directed graph $\Gamma$ underlying the boustrophedon transform.}
\end{figure}
Let $\pi (n,k,i)$ denote the number of paths in $\Gamma$ from the node
labeled $T_{i, 0}$ to the node labeled $T_{n,k}$.
It follows from the rule for constructing the triangle that the numbers
$T_{n,k}$ are given by
\setcounter{equation}{6}
\beql{Eq11}
T_{n,k} = \sum_{i=0}^n \pi (n,k,i) a_i ~.
\eeq

From Section~1 we know that the boustrophedon transform of the sequence
$1,0,0,0, \dd$ is $E_0$, $E_1$, $E_2$, $E_3$, $\dd$, and so
(from \eqn{Eq11})
\beql{Eq12}
E_n = \pi (n,n,0) ~~~~ (n \ge 0 ) ~.
\eeq
We will give a direct proof of this (although of course it is known result,
cf. \cite{Arn92}), in order to establish a bijection
between paths in $\Gamma$ and up-down permutations.
\begin{prop}
\label{pr1}
$\pi (n,n,0)$ is equal to $E_n$,
(\htmladdnormallink{A1111}{http://www.research.att.com/cgi-bin/access.cgi/as/njas/sequences/eisA.cgi?Anum=000111}),
the number of down-up permutations
of $\{1,2, \dd , n \}$.
\end{prop}
\paragraph{Proof.}
Let $P$ be a path in $\Gamma$ from the top node to the node labeled
$T_{n,n}$.
(Fig.~2 shows an example for $n=5$.)
\begin{figure}[htb]
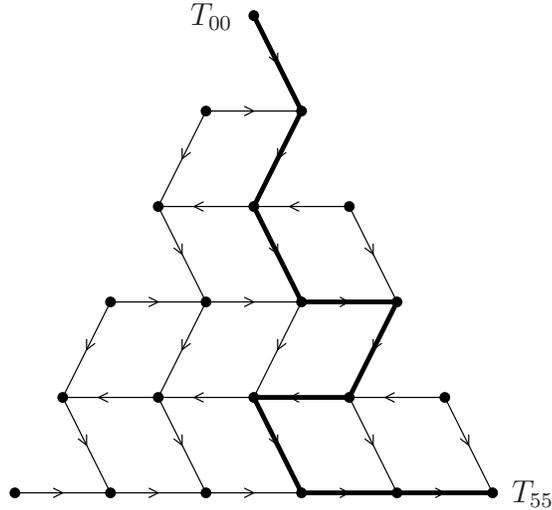

\begin{center}
\input bo2.tex
\end{center}
\caption{A path from $T_{0,0}$ to $T_{5,5}$.}
\end{figure}
Let $T_{i,f(i)}$ be the label of the node where $P$ arrives at level $i$
$(1 \le f(i) \le i \le n )$.
We construct a box diagram to represent $P$ by the following
procedure (see Fig.~3).
The bottom row contains $n$ boxes labeled $1, \dd, n$ from left to
right (if $n$ is even) or from right to left (if $n$ is odd).
The box labeled $f(n)$ is starred.
We now repeatedly place a row of boxes above the empty boxes,
putting a star in the $f(i)$-th box, always counting from the left if
$i$ is even or from the right if $i$ is odd,
for $i=n-1, n-2, \dd, 1$.
\begin{figure}[htb]
\centerline{\psfig{file=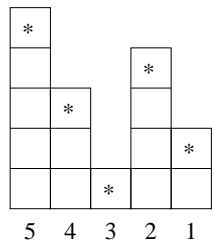,width=1.1in}}
\caption{Box diagram corresponding to path $P$ in Fig.~2.}
\end{figure}

We convert the box diagram into a permutation of $\{1, \dd, n\}$ by reading the
rows from the bottom up and recording the number at the foot of the column
containing the star.
(The permutation corresponding to the above example is
$(3,1,4,2,5)$.)
We omit the easy verification that this process defines a bijection between
paths and down-up permutations. \hfill $\bsq$
\begin{prop}
\label{pr2}
$$\pi (n,n,k) = {\binom{n}{k}} E_{n-k} ,~~~\mbox{for}~~
0 \le k \le n ~.$$
\end{prop}
\paragraph{Sketch of proof.}
Consider a path from the node labeled $T_{k,0}$ to the node labeled $T_{n,n}$,
such as the path from $T_{4,0}$ to $T_{9,9}$ shown in Fig.~4.
\begin{figure}[htb]
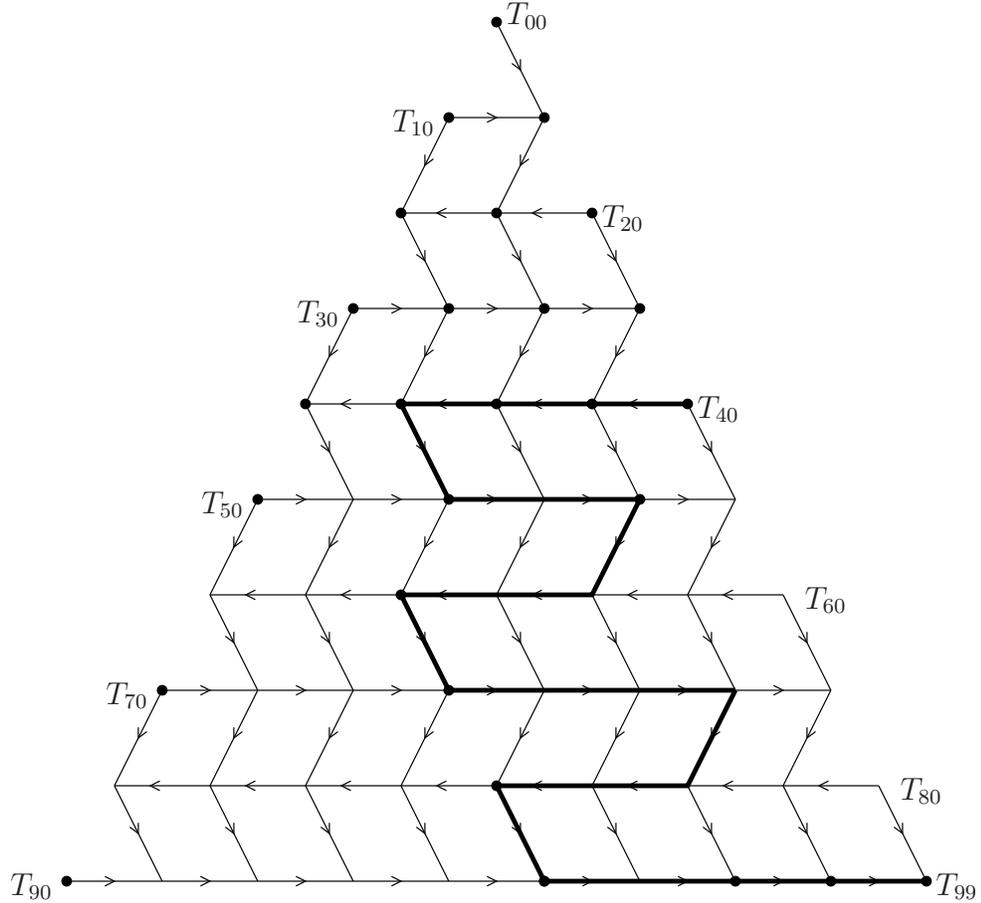

\begin{center}
\input bo4.tex
\end{center}
\caption{A path from $T_{4,0}$ to $T_{9,9}$.}
\end{figure}
The procedure used in the proof of Proposition~\ref{pr1} converts this into
a box diagram,
which for this example is shown in Fig.~5.
\begin{figure}[htb]
\centerline{\psfig{file=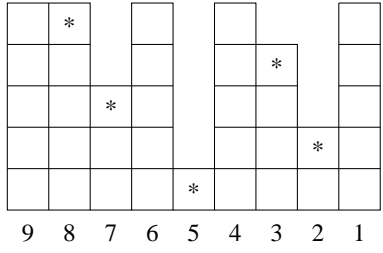,width=2in}}
\caption{Box diagram for path shown in Fig.~4.}
\end{figure}
The columns that do not contain stars identify one of the ${\binom{n}{k}}$
$k$-subsets of $\{1, \dd, n \}$, while the starred columns
themselves form a box diagram (in this case it is that shown in
Fig.~3) that identifies a down-up permutation of $\{1, \dd, n-k \}$. \hfill $\bsq$

From Proposition~\ref{pr2} and \eqn{Eq11} we obtain
$$b_n = T_{n,n} = \sum_{k=0}^n {\binom{n}{k}} E_{n-k} a_k ~,$$
which establishes \eqn{Eq6}.
Equations \eqn{Eq8} and \eqn{Eq7} now follow immediately.
This completes the proof of the theorem. \hfill $\bsq$
\paragraph{Remark.}
With only a little more effort we can determine all the ``boustrophedon numbers''
$\pi (n,k,i)$.
Note that $\pi(n,0,i) =0$ for $n \ge 1$, $0 \le i \le n-1$, and $\pi (n,0,n) =1$.
\begin{prop}
\label{pr3}
For $n \ge 1$, $0 \le k \le n-1$,
$$\pi (n,k,0) = E_{n,k} = \sum_{r=0}^{[(k-1)/2]} (-1)^r
{\binom{k}{2r+1}} E_{n-2r-1} ~.$$
\end{prop}
\paragraph{Proof.}
$\pi(n,k,0) = E_{n,k}$ follows from \eqn{Eq11} and the definition of
$E_{n,k}$ (see \eqn{Eq2}), and the formula for $E_{n,k}$ is given in \cite{Ent66}. \hfill $\bsq$
\paragraph{Remark.}
If the path is extended to reach the node labeled $T_{n+1, n+1}$, the
corresponding box diagram has the same format as those arising in
Proposition~\ref{pr1}, except that the star in the last row is constrained
to appear in the box labeled $k+1$.
\begin{prop}
\label{pr4}
For $n \ge 2$, $0 < k < n$, $0 < i \le n$,
\beql{Eq20}
\pi (n,k,i) = \sum_{s=0}^{\min \{k, n-i\}}
{\binom{k}{s}} {\binom{n-k}{n-i-s}}
\pi (n-i,s,0) ~.
\eeq
\end{prop}
\paragraph{Sketch of proof.}
Consider a path $P$ from $T_{i,0}$ to $T_{n,k}$, and complete it to a path $Q$ from $T_{i,0}$ to $T_{n+1, n+1}$ by extending $P$ by a downward sloping edge and a series of
horizontal edges, as illustrated in Fig.~6.
\begin{figure}[htb]
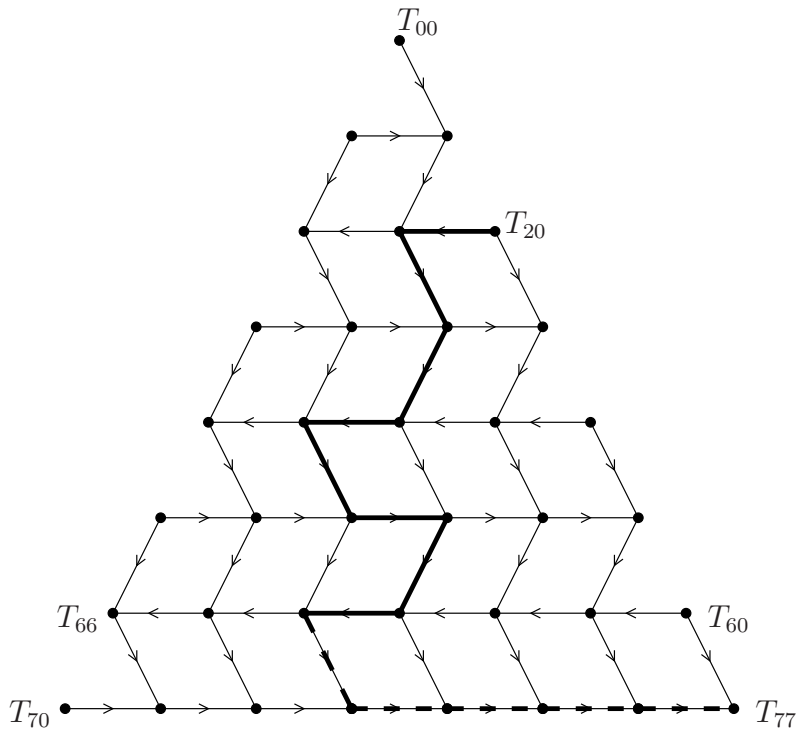

\begin{center}
\input bo6.tex
\end{center}
\caption{Path from $T_{2,0}$ to $T_{6,4}$ (solid line)
and its continuation to $T_{7,7}$ (broken line).}
\end{figure}
We form the box diagram for $Q$, as in Proposition~\ref{pr2} (see Fig.~7a).
After deleting all the unstarred columns we obtain the box
diagram for a path of type $\pi (n-i,s,0)$,
for some $s$ (Fig.~7b).
\begin{figure}[H]
\centerline{\psfig{file=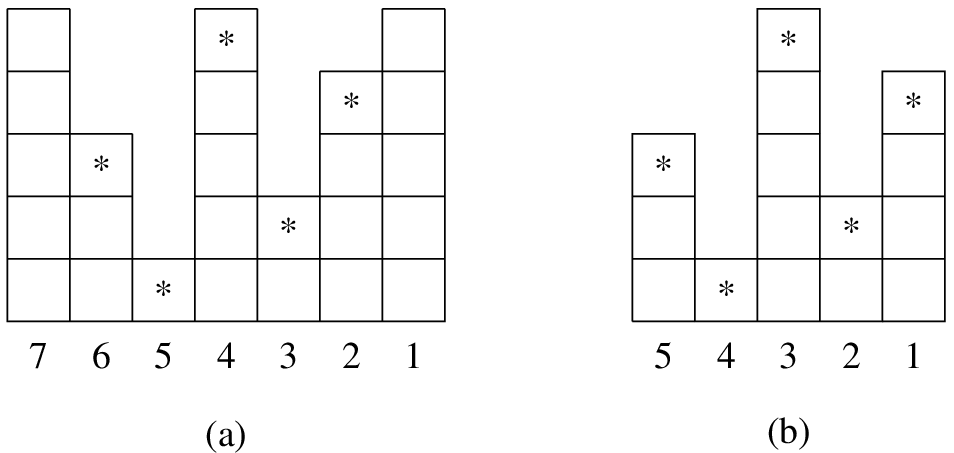,width=3in}}
\caption{(a)~Box diagram corresponding to path shown in Fig.~6.
(b)~Reduced box diagram, representing path of type $\pi (4,3,0)$.}
\end{figure}

In the box diagram for $Q$ itself, the star in the last row divides
the remaining stars into two sets of sizes $s$ (to the right) and $n-i-s$ (to the left),
and the binomial coefficients in \eqn{Eq20} count the ways in which the corresponding columns can be selected. \hfill $\bsq$

Propositions~\ref{pr1}--\ref{pr4} together express all the boustrophedon numbers
in terms of the $E_n$'s,
and via \eqn{Eq11} give an explicit formula for every entry in the triangle
\eqn{Eq9}.
\section{Combinatorial interpretations and examples}
\hsp
Equation~\eqn{Eq6} yields
many possible combinatorial interpretations for the numbers $b_n$.
For example, if $a_n$ is the number of arrangements of $n$ labeled
objects so that they have some property $Q$, then $b_n$ is the number of ways
of dividing $n$ objects into two groups so that the first group has property $Q$
and the second forms a down-up sequence.
Since $E_n$ is also the number of ordered binary trees on $n$ nodes
(cf. \cite{Pou82}, \cite{KPP94}), other interpretations for the $b_n$ can be given in terms of graphs.
\begin{exam}
\label{ex1}
{\rm
We can see now that \eqn{Eq15} has e.g.f. $e^x ( \mbox{sec} \, x + \tan x )$, and that the $n$-th term of this sequence gives the
number of ways that we can form a down-up sequence of some length $\ell \ge 0$
from $\{1, \dd, n \}$.
E.g. for $n=3$ there are 9 possibilities:
$\phi , 1,2,3,21,31,32,213,312$.
}
\end{exam}
\begin{exam}
\label{ex2}
{\rm
The boustrophedon transform of the Bell numbers (cf. \cite{SP95}, Fig.~M4981;
sequence \htmladdnormallink{A110}{http://www.research.att.com/cgi-bin/access.cgi/as/njas/sequences/eisA.cgi?Anum=000110});
produces the sequence $1,2,5,16,60,258, \dd$,
(\htmladdnormallink{A764}{http://www.research.att.com/cgi-bin/access.cgi/as/njas/sequences/eisA.cgi?Anum=000764})
whose $n$-th term gives the
number of ways to take blocks labeled $1, \dd, n$ and to partition some of them into heaps and to arrange the rest so they form a down-up sequence.
}
\end{exam}
\begin{exam}
\label{ex3}
{\rm
The boustrophedon transform of the $E_n$ sequence shifted one place to the left is the same sequence
shifted two places to the left:
$$
\begin{array}{ccccccccc}
~ & ~ & ~ & ~ & 1 \\
~ & ~ & ~ & 1 & \to & 2 \\
~ & ~ & 5 & \leftarrow & 4 & \leftarrow & 2 \\
~ & 5 & \to & 10 & \to & 14 & \to & 16 \\
61 & \leftarrow & 56 & \leftarrow & 46 & \leftarrow & 32 & \leftarrow & 16 \\
~ & ~ & ~ & ~ & \cdots
\end{array}
$$
In view of Theorem~\ref{th1}, this means the e.g.f. $\sE (x)$ satisfies
$$\sE (x) \sE ' (x) = \sE '' (x) ~.$$
The initial conditions $E_0 = E_1 = 1$ then give $\sE (x) = {\rm sec} \, x + \tan x$ as the solution.
}
\end{exam}
\begin{exam}
\label{ex4}
{\rm
The sequence $1, 0, 1, 1,2,6,17,62,259,1230, \dd$
(\htmladdnormallink{A661}{http://www.research.att.com/cgi-bin/access.cgi/as/njas/sequences/eisA.cgi?Anum=000661}) is the
lexicographically earliest sequence that begins with 1 and shifts
{\em two} places left under the boustrophedon transform.
(Examples~\ref{ex3} and \ref{ex4} are both eigen-sequences for this
transform, in the notation of \cite{BS95}.)
We do not know of any combinatorial interpretation for these numbers.
}
\end{exam}
\begin{exam}
\label{ex5}
{\bf The double-ox transform.}
{\rm
Generalizing some examples of Arnold (\cite{Arn92}, see also \cite{Dum95}), we consider two oxen
plowing separate fields with a messenger that takes the output at the end of one row and rushes it to be used by the other
ox as input to the next row.
For example, if the initial sequence (shown in italics in Fig.~8) is
$1,1,1, \dd$, this produces the output sequence (shown in bold)
$1,3,9,35,177,1123, \dd$
(\htmladdnormallink{A834}{http://www.research.att.com/cgi-bin/access.cgi/as/njas/sequences/eisA.cgi?Anum=000834}).
\begin{figure}[htb]
$$
\begin{array}{ccccccccc}
~ & ~ & ~ & ~ & {\mbox{\em 1}} \\
~ & ~ & ~ & {\mathbf 3} & ~ & 2 \\
~ & ~ & {\mbox{\em 1}} & ~ & 4 & ~ & 6 \\
~ & {\mathbf 35} & ~ & 34 & ~ & 30 & ~ & 24 \\
{\mbox{\em 1}} & ~ & 36 & ~ & 70 & ~ & 100 & ~ & 124 \\
~ & ~ & ~ & ~ & \cdots
\end{array}~~~~~~~
\begin{array}{ccccccccc}
~ & ~ & ~ & ~ & {\mathbf 1} \\
~ & ~ & ~ & 2 & ~ & {\mbox{\em 1}} \\
~ & ~ & 6 & ~ & 8 & ~ & {\mathbf 9} \\
~ & 24 & ~ & 18 & ~ & 10 & ~ & {\mbox{\em 1}} \\
124 & ~ & 148 & ~ & 166 & ~ & 176 & ~ & {\mathbf 177} \\
~ & ~ & ~ & ~ & \cdots
\end{array}
$$
\caption{The double-ox transform of $1,1,1, \dd$ is
$1,3,9,35,177, \dd$.}
\end{figure}
}
\end{exam}

Less colorfully, let $\ba = a_0, a_1, \dd$ be the initial
sequence,
$\bm = m_0 , m_1, \dd$ the middle (or messenger) sequence, and
$\bb = b_0 , b_1, \dd$ the transformed sequence.
We define two triangles of numbers $\{ L_{n,k} \}$ and $\{R_{n,k} \}$, with $0 \le k \le n$, by
\begin{eqnarray*}
L_{2i,0} & = & a_{2i} ,~~~
R_{2i+1,0} ~=~ a_{2i+1} , \\ [+.05in]
L_{2i,2i} & = & R_{2i,0} ~=~ m_{2i} , ~~~
L_{2i+1,0} ~=~ R_{2i+1 ,2i+1} ~=~ m_{2i+1} \\ [+.05in]
L_{2i+1, 2i+1} & = & b_{2i+1} , ~~~R_{2i,2i} ~=~ b_{2i} , \\
\noalign{and}
L_{n+1, k+1} & = & L_{n+1, k} + L_{n,n-k} , ~~~
R_{n+1,k+1} ~=~ R_{n+1,k} + R_{n,n-k} ~.
\end{eqnarray*}

We were happy to find that Theorem~\ref{th1} leads to an equally simple
description of this transformation.
The proof is left to the reader.
\begin{theorem}
\label{th2}
The e.g.f.'s of $\ba$, $\bm$ and $\bb$ are related by
\begin{eqnarray*}
\sM (x) & = & \frac{1}{\cos x - \sin x} ~ \sA (x) ~, \\ [+.1in]
\sB (x) & = & \frac{\cos x + \sin x}{\cos x - \sin x} ~
\sA (x) ~.
\end{eqnarray*}
\end{theorem}

\subsection*{Acknowledgements}
\hsp
We thank Richard Guy for suggesting the problem, and J.~H. Conway,
I.~M. Pak and S. Plouffe for helpful discussions.
A referee suggests that since the boustrophedon and double-ox
transforms are associated with the root systems $A_n$ and $C_n$,
respectively (cf. \cite{Arn92}), it would be worthwhile studying analogues
of these transforms for other root systems.

\subsection*{Postscript added June 10 2002}
\hsp
(1)  David Callan points out that the proof of Proposition~\ref{pr1} is not quite right.
He says:

\begin{quote}
The bijection from paths to down-up permutations is not 
quite right. (It is possible for $f(n)$ to equal 1 and then the box diagram
procedure as described would give a permutation beginning with
1---definitely not down-up.)
One way to correct it would be the following.
Label the $n$ boxes in the bottom row $1,2,3,\ldots,n$ from left to right in all
cases (whether $n$ is even or odd). Then place stars as described 
in the present proof,
except counting from the right for the first (bottom-most) star, then
alternately left and right for subsequent stars.
\end{quote}

We agree with this comment, and thank him for the correction.

\hsp

(2) Mike Atkinson drew our attention to two earlier
references that we could have cited:

M.D. Atkinson: Zigzag permutations and comparisons of adjacent elements,
{\em Information Processing Letters}, {\bf 21} (1985), 187--189.

M.D. Atkinson: Partial orders and comparison problems, in {\em Sixteenth
Southeastern Conference on Combinatorics, Graph Theory and Computing, (Boca
Raton, February 1985)}, Congressus Numerantium {\bf 47}, 77--88.

\end{document}